\documentclass[review,number]{elsarticle}
\usepackage[utf8]{inputenc}
\usepackage[T1]{fontenc}
\usepackage{lmodern}% use this package to use modern version of cm fonts (useable with T1 encoding)
\usepackage[english]{babel}
\usepackage{amsfonts,amsmath,amssymb}
\usepackage{graphicx}
\usepackage{amsthm}
\usepackage{xcolor}
\usepackage{url}

\usepackage{eso-pic}
\usepackage{pict2e}

\newtheorem{theorem}{Theorem}[section]
\newtheorem{lemma}[theorem]{Lemma}

\DeclareMathOperator{\id}{Id}
\DeclareMathOperator{\prox}{prox}
\DeclareMathOperator{\grad}{grad}

\DeclareMathOperator{\TV}{TV}
\DeclareMathOperator{\ind}{ind}
\DeclareMathOperator*{\argmin}{arg\,min}
\newcommand{\realr}{\mathbb{R}}
\newcommand{\naturaln}{\mathbb{N}}
\newcommand{\realrplusinfty}{\ensuremath{{\mathbb{R}_{+\infty}}}}
\newcommand{\stackvec}[2]{\left(\begin{array}{c} #1\\ #2 \end{array}\right)}

\newcommand{\amendment}[1]{{\color{black} #1}}% final (all black)

\AddToShipoutPicture*{%

  \AtTextUpperLeft{\parbox[b]{15cm}{©2024. This manuscript version is made available under the CC-BY-NC-ND 4.0 license.\\
\url{https://creativecommons.org/licenses/by-nc-nd/4.0/}\\[3cm]}}%
}

\begin{document}

\title{Convergence analysis of a primal-dual optimization-by-continuation algorithm}

\author[1]{Ignace Loris\corref{cor1}}
\ead{Ignace.Loris@ulb.be}
%\ead[url]{ignace.loris.web.ulb.be}

\author[2]{Simone Rebegoldi}
\ead{simone.rebegoldi@unimore.it}
\cortext[cor1]{Corresponding author}

%\fntext[fn1]{This is the first author footnote.}
%\fntext[fn2]{Another author footnote, this is a very long footnote and it should be a really long footnote. But this footnote is not yet sufficiently long enough to make two lines of footnote text.}
%\fntext[fn3]{Yet another author footnote.}

\affiliation[1]{organization={Mathematics Department, Université libre de Bruxellles},
addressline={Boulevard du Triomphe},
postcode={1050},
city={Brussels},
country={Belgium}}

\affiliation[2]{organization={Dipartimento di Scienze Fisiche, Informatiche e Matematiche, Università degli
Studi di Modena e Reggio Emilia},
addressline={Via Campi 213/b},
city={Modena},
postcode={41125},
country={Italy}}

\begin{abstract}
We present a numerical iterative  optimization algorithm for the minimization of a cost function consisting of a linear combination of three convex terms, one of which is differentiable, a second one is prox-simple and the third one is the composition of a linear map and a prox-simple function. The algorithm's special feature lies in its ability to approximate, in a single iteration run, the minimizers of the cost function for many different values of the parameters determining the relative weight of the three terms in the cost function. A proof of convergence of the algorithm, based on an inexact variable metric approach, is also provided. As a special case, one recovers a generalization of the primal-dual algorithm of Chambolle and Pock, and also of the proximal-gradient algorithm. Finally, we show how it is related to a primal-dual iterative algorithm based on inexact proximal evaluations of the non-smooth terms of the cost function.
\end{abstract}

%%Research highlights
%\begin{highlights}
%\item Research highlight 1
%\item Research highlight 2
%\end{highlights}

\begin{keyword}
primal-dual algorithm \sep proximal algorithm \sep inverse problem \sep iterative algorithm \MSC[2020]{90C06}\sep\MSC[2020]{90C25}\sep\MSC[2020]{90C59}\sep\MSC{90C90}\sep\MSC[2020]{49M29}\sep\MSC[2020]{65K10}
\end{keyword}
\maketitle

\section{Introduction}\label{sec-intro}

We study the optimization problem
\begin{equation}\label{eq:problem}
\argmin_{u\in\realr^{d}} F_{\lambda,\mu}(u)\stackrel{\text{def.}}{=} f(u)+\lambda g(u)+\mu h(Au)
\end{equation}
with the following assumptions
\begin{itemize}
\item $f:\realr^d\to\realr$ is differentiable with Lipschitz continuous gradient (constant $L$);
\item $g:\realr^{d}\to\realrplusinfty$ and $h:\realr^{d'}\to\realrplusinfty$ are proper, convex and closed functions (here we use $\realrplusinfty=\realr\cup\{+\infty\}$);
\item $A:\realr^d\to\realr^{d'}$ is a linear map;
\item positive parameters $\lambda,\mu>0$;
\item there exists at least one minimizer $\hat u_{\lambda,\mu}$ of $F_{\lambda,\mu}$.
\end{itemize}
Furthermore, and from a practical side, it is also assumed that the gradient of $f$ and the proximal operator (see expression~(\ref{eq:prox-def})) of $g$ and $h$ can be evaluated. It is \emph{not} assumed that the proximal operator of $h\circ A$ is available.

We will be interested in the case when the number of unknowns $d$ is large. Such optimization problems are often encountered in imaging and other inverse problems \cite{Chambolle2016,Bottou2018,Bach2012}. In that context, the first term in the cost function (\ref{eq:problem}) typically corresponds to a data fidelity term (e.g. as in least squares) and the two remaining terms correspond to penalties that need to be added to discard `unphysical' solutions arising from the ill-posedness of the underlying inverse problem \cite{Bertero2021}.

In this context the values of the parameters $\lambda$ and $\mu$, which regulate the relative importance of each term in the cost function (\ref{eq:problem}), are not necessarily known in advance. In other words, one is interested in computing the minimizers $\hat u_{\lambda,\mu}$ for many different values of those parameters.

It can be shown that the minimizers $\hat u_{\lambda,\mu}$ trace out a frontier (a curve or surface) when plotting the fidelity term $\sigma=f(\hat u_{\lambda,\tau})$ in terms of the penalty terms $\tau_1=g(\hat u_{\lambda,\mu})$ and $\tau_2=h(A\hat u_{\lambda,\mu})$ for different values of $\lambda$ and $\mu$. This `Pareto frontier' indicates where no further improvement can be made to the fidelity term without making at least one of the two penalty terms larger. In other words, it defines the boundary (in the space of $(\tau_1,\tau_2;\sigma)\in\realr^3$) between the region above it, which can be accessed by points of the form $(g(u),h(Au);f(u))$ for $u\in\realr^d$, and the region below it which is inaccessible by points of the form $(g(u),h(Au);f(u))$. 
Section~\ref{sec-trade-off-surface} and Theorem~\ref{prop:Pareto} review some properties (convexity and subddiferentiability) of the Pareto frontier in the context of problem (\ref{eq:problem}).

The main contribution of this paper is the introduction of an iterative algorithm that aims at computing (approximately) the minimizers $\hat u_{\lambda,\mu}$ for many different values of the parameters $\lambda$ and $\mu$. 
%In this paper, we propose an iterative first order algorithm that aims to do just that. 
The proposed algorithm is a first order primal-dual method related to the one of \cite{Condat2013}, but adapted in such a way that its intermediate iterates can be interpreted as good approximations of solutions of problem (\ref{eq:problem}) for different values of $\lambda$ and $\mu$. 
In particular, two converging sequences $\lambda_n\to\lambda$ and $\mu_n\to\mu$ are introduced, and the proposed iterative algorithm uses these variable penalty parameters at every iteration step, instead of the fixed parameters $\lambda, \mu$. This gives rise to the algorithm
\begin{equation}\label{eq:algorithm}
\left\{
\begin{array}{lcl}
u_{n+1} & = & \prox_{\alpha\lambda_n g}\left(u_n-\alpha \grad f(u_n)-\alpha\mu_n A^T  v_n\right)\\
v_{n+1} & = & \prox_{(\beta/\mu_n)\, h^\ast}\left(v_n+(\beta/\mu_n) A (2u_{n+1}-u_n)\right)
\end{array}
\right.\quad n= 0,1,\ldots
\end{equation}
with any initial point $(u_0,v_0)\in\realr^{d+d'}$ and where $v_n$ are auxiliary (dual) variables. Under suitable conditions on the step length parameters $\alpha$ and $\beta$, we show that this algorithm still converges to a minimizer of problem (\ref{eq:problem}).

The main benefit of algorithm (\ref{eq:algorithm}) lies in the additional freedom of choosing the sequences $\lambda_n$ and $\mu_n$. If one can choose $\lambda_0$ and $\mu_0$ large and $u_0\approx\hat u_{\lambda_0,\mu_0}$, and at the same time decrease $\lambda_n$ and $\mu_n$ slowly, one may expect that the iterates $u_n$ are all good approximations of each $\hat u_{\lambda_n,\mu_n}$. In this way, one can approximate, in a single iteration run, the minimizers of the cost function (\ref{eq:problem}) for many values of the parameters $\lambda$ and $\mu$. Such an optimization-by-continuation type algorithm has already been introduced before for a simplified problem ($\mu=0$), and without proof of convergence in \cite{Hale2008}. Here the main difference is the inclusion of a second penalty term with non prox-simple penalty function, the primal-dual nature of the algorithm and also the inclusion of a proof of convergence.  As a special case, we also recover a minimization-by-continuation version of the well-known Chambolle-Pock algorithm \cite{Chambolle2011} popular in imaging.

From a technical perspective, the proof technique is somewhat different from the method used in \cite{Fest2023} (which treats the special case $\mu=0$ of problem (\ref{eq:problem})), as it relies here on the introduction of a variable metric. It is also different (more explicit) from the proof of convergence given in \cite{Condat2013} (for the case $\lambda_n=\lambda$ and $\mu_n=\mu$ fixed) which is based on a reduction to the Krasnosel'ski\v{\i}-Mann iteration theorem for averaged operators \cite{Dong2022,Ryu2022}.

We also show how the proposed algorithm (\ref{eq:algorithm}) could be interpreted as the primal--dual method with inexact proximal evaluations proposed in \cite{Rasch2020}. 
%However, stricter conditions need to be imposed on the cost function in order to rely on the convergence results of \cite{Rasch2020}.

Section~\ref{sec-trade-off-surface} discusses the properties of the so-called Pareto frontier attached to problem (\ref{eq:problem}). Section~\ref{sec-pd-alg} contains the derivation and proof of convergence of the primal-dual optimization-by-continuation algorithm (\ref{eq:algorithm})  and the relation to inexact proximal evaluations is discussed in Section~\ref{sec:inexact}. Section~\ref{sec:numerics} contains two numerical illustrations of the proposed algorithm in the area of image reconstruction. \amendment{The Python code used for the numerical experiments is available online \cite{Loris2024b}.}

%origin of prox grad: \cite{Fukushima1981} see also \cite{Kanzow2021} for this remark.
%
%this can be reinterpreted as a single step of the usual algorithm for $\lambda_n,\mu_n$

%, introduce some freedom in existing method, goal is to choose a better path (somehow?), paper proves that convergence is still ok .

Finally, we briefly introduce some definitions and notations used throughout the paper. We work with the set $\Gamma_0(\realr^d)$ of proper, convex and closed function defined on $\realr^d$. We shall make use of the Legendre-Fenchel conjugate of a function $f\in\Gamma_0(\realr^d)$ defined by
\begin{equation}
f^\ast:\realr^d\to\realrplusinfty: \xi\to f^\ast(\xi)=\sup_{u\in\realr^d}\langle\xi,u\rangle-f(u)
\end{equation}
and also of the sub-differential $\partial f$ defined by:
\begin{equation}
\partial f(a)=\{\xi\in\realr^d\text{ such that } f(u)\geq f(a)+\langle\xi,u\rangle\quad \forall u\in\realr^d\}.
\end{equation}
On a technical level, we assume in this paper that $g$, $A$ and $h$ appearing in problem (\ref{eq:problem}) are such that $\partial (\lambda g+\mu h\circ A)=\lambda \partial g+\mu A^T \partial h\circ A$ (see e.g. \cite[Th. 23.8 and 23.9]{Rockafellar1970} or \cite[Lemma~2.2]{Combettes2005}).
Furthermore we shall use the proximal operator \cite{Moreau1965,Combettes2005} of a function $f\in\Gamma(\realr^d)$ which is defined as:
\begin{equation}\label{eq:prox-def}
\prox_f:\realr^d\to\realr^d: a\to\prox_f(a)=\arg\min_{u\in\realr^d}\frac{1}{2}\|u-a\|_2^2+f(u).
\end{equation}
This operator is related to the proximal operator of the Fenchel conjugate function $f^\ast$ by the identity $\prox_{f}(u)+\prox_{f^\ast}(u)=u$ ($\forall u\in\realr^d$) or more generally:
\begin{equation}
u=\prox_{\alpha f}(u)+\alpha\prox_{\alpha^{-1}f^\ast}(\alpha^{-1}u).
\end{equation}
It can be shown that $\prox_{\alpha f}(x)$ is continuous with respect to $(x,\alpha)$ for any $x\in\realr^d$ and $\alpha>0$.
These proximal operators can also be used to characterize the sub-gradient $\xi$ of $f$ in $u$ as \cite{Moreau1965}:
\begin{equation}\label{eq:subdiffcharac}
\begin{array}{clcl}
&\xi\in\partial f(u) &\Leftrightarrow& u=\prox_{\alpha f}(u+\alpha \xi) \text{ for any } \alpha>0 \\[3mm]
\Leftrightarrow& u\in\partial f^\ast(\xi) &\Leftrightarrow & 
\xi=\prox_{\beta f^\ast}(\beta u+\xi) \text{ for any }  \beta>0. %\\[3mm]
%\Leftrightarrow& f(u)+f^\ast(\xi)=\langle \xi,u\rangle.
\end{array}
\end{equation}
We have written the standard euclidean norm as $\|x\|_2$ and norm with respect to a positive definite matrix $M$ as $\|x\|_M$. The spectral norm of a matrix will be written $\|M\|$ (largest singular value of $M$).

\section{Properties of the Pareto frontier}\label{sec-trade-off-surface}

In order to discuss the properties of the Pareto frontier in the setting of problem (\ref{eq:problem}), we follow and expand upon the analysis of \cite{Berg2008,Berg2011} who treat a subcase of the special case $\mu=0$. First we introduce the closely related constrained problem
\begin{equation}\label{eq:constrainedproblem}
\min_{u\in\realr^d} f(u) \quad\text{with}\quad g(u)\leq \tau_1,\ h(Au)\leq \tau_2
\end{equation}
and denote $\tau = (\tau_1,\tau_2)$ and $\tilde u_{\tau}$ the minimizer of problem (\ref{eq:constrainedproblem}) which we assume to exist when the feasible set
\begin{displaymath}
S(\tau)\stackrel{\text{def.}}{=}\{u\in\realr^d, g(u)\leq \tau_1, h(Au)\leq \tau_2\}
\end{displaymath} 
is non-empty.

The properties of the Pareto frontier of problem (\ref{eq:problem}) can now be most easily derived in terms of the \emph{value function}:
\begin{equation}\label{eq:valuefunction}
\varphi:\realr^2\to\realrplusinfty:\tau\to\varphi(\tau)\stackrel{\text{def.}}{=}\inf\left\{ f(u) \quad\text{with}\quad g(u)\leq \tau_1,\ h(Au)\leq \tau_2\right\}
\end{equation}
of the constrained problem (\ref{eq:constrainedproblem}). In fact, we define the Pareto frontier (or trade-off frontier) of problem (\ref{eq:problem}) as the graph of value function (\ref{eq:valuefunction}), i.e. $\{ (\tau_1,\tau_2,\varphi(\tau_1,\tau_2)),\tau\in\realr^2\}$. The properties listed in Theorem ~\ref{prop:Pareto} below justify this definition. Their proof makes use of the Lagrange dual problem \cite{Boyd2004} of the constrained problem (\ref{eq:constrainedproblem}) which is defined by
\begin{equation}\label{eq:dual}
\sup_{\lambda,\mu\geq 0} \tilde f(\lambda,\mu)\quad\text{with}\quad \tilde f(\lambda,\mu)\stackrel{\text{def.}}{=} \inf_{u\in\realr^d}f(u)+\lambda (g(u)-\tau_1)+\mu (h(Au)-\tau_2).
\end{equation}

\begin{theorem}\label{prop:Pareto} If $f,g,h:\realr^d\to\realrplusinfty$ are convex functions, one has:
\begin{enumerate}
\item The value function $\varphi:\realr^2\to\realrplusinfty$ is non-increasing; 
\item The value function $\varphi:\realr^2\to\realrplusinfty$ is convex;
\item The region below the graph $(\tau,\varphi(\tau))_{\tau\in\realr^2}$ cannot be reached by any point of the form $((g(u),h(Au)),f(u))$ with $u\in\realr^d$;

\item\label{item:duality} If $\lambda,\mu> 0$ and $\hat u_{\lambda,\mu}\in\arg\min_{u \in \realr^d} f(u)+\lambda g(u)+\mu h(Au)$ then $\hat u_{\lambda,\mu}\in \arg\min_{u\in S(\tau)} f(u)$ for $\tau=(g(\hat u_{\lambda,\mu}),h(A\hat u_{\lambda,\mu}))$ and strong duality holds: $f(\hat u_{\lambda,\mu})=\tilde f(\lambda,\mu)$.

\item If $\lambda,\mu> 0$  the point $((g(\hat u_{\lambda,\mu}),h(A\hat u_{\lambda,\mu})),f(\hat u_{\lambda,\mu}))$  is on the Pareto frontier;

\item If $\lambda,\mu> 0$ the subdifferential of the value function at the point $\tau=(g(\hat u_{\lambda,\mu}),h(A\hat u_{\lambda,\mu}))$ contains $-(\lambda,\mu)$: $-(\lambda,\mu)\in\partial\varphi(\tau)$;

\item If $-(\lambda,\mu)\in\partial \varphi(\tau)$, then $\tilde u_{\tau}\in\arg\min F_{\lambda,\mu}(u)$.

%\item The function $\varphi$ is differentiable on the interior of its domain???? needs more!! counterexample??
%the domain of $\varphi$ and the derivative in the point $(g(\hat u_{\lambda,\mu}),f(\hat u_{\lambda,\mu}))$ is $-(\lambda,\mu)$.
\end{enumerate}
\end{theorem}
\begin{proof}
1) If $\tau\geq \tilde\tau$ (componentwise) the feasible set $S(\tau)$ contains the feasible set $S(\tilde\tau)$. Therefore $\inf\left\{ f(u),\ u\in S(\tau)\right\}\leq \inf\left\{ f(u), \ u\in S(\tilde\tau)\right\}$.

2) In order to show convexity of the value function, we follow \cite[p.~50]{Ekeland1999}. By definition of $\varphi$ one has:
\begin{displaymath}
\begin{array}{l}
\forall a \text{ such that } \varphi(\tau)<a\quad \exists u\in S(\tau) \text{ such that } \varphi(\tau)\leq f(u)<a\\[3mm]
\forall \tilde a \text{ such that } \varphi(\tilde\tau)<\tilde a\quad \exists \tilde u\in S(\tilde\tau) \text{ such that } \varphi(\tilde\tau)\leq f(\tilde u)<\tilde a.
\end{array}
\end{displaymath}
Choosing any $\alpha\in[0,1]$, one finds $g(\alpha u+(1-\alpha) \tilde u)\leq \alpha g(u)+(1-\alpha) g(\tilde u)\leq \alpha \tau+(1-\alpha) \tilde \tau$ and likewise $h(A (\alpha u+(1-\alpha) \tilde u))\leq \alpha \tau+(1-\alpha) \tilde \tau$. Therefore $\alpha u+(1-\alpha) \tilde u\in S(\alpha \tau+(1-\alpha) \tilde \tau)$ by convexity of $g$ and $h$.
Now we can write:
\begin{displaymath}
\begin{array}{lcl}
 \displaystyle \varphi(\alpha \tau+(1-\alpha)\tilde\tau) & \stackrel{(\ref{eq:valuefunction})}{=} & \displaystyle 
\inf\left\{ f(v) \quad\forall v\in S(\alpha \tau+(1-\alpha)\tilde\tau)\right\}
\\[3mm]
&\leq& f(\alpha u+(1-\alpha)\tilde u)\\[3mm]
&\stackrel{f\text{ convex}}{\leq}& \alpha  f(u)+(1-\alpha)f(\tilde u)\\[3mm]
&<& \alpha  a+(1-\alpha)\tilde a.
\end{array}
\end{displaymath}
By decreasing $a$ to $\varphi(\tau)$ and $\tilde a$ to $\varphi(\tilde \tau)$ one obtains
\begin{displaymath}
\varphi(\alpha \tau+(1-\alpha)\tilde\tau)\leq \alpha  \varphi(\tau)+(1-\alpha)\varphi(\tilde\tau).
\end{displaymath}

3) If a point $u\in\realr^d$ with $(g(u),h(Au))=\tau$ and $f(u)<\varphi(\tau)$ were to exist, then this would be in contradiction with the definition (\ref{eq:valuefunction}) of the value function $\varphi$.

4) Let $\lambda,\mu\geq 0$, set $\hat u_{\lambda,\mu}\in \arg\min_{u\in\realr^d} f(u)+\lambda g(u)+\mu h(Au)$ and $\tau=(g(\hat u_{\lambda,\mu}),h(A \hat u_{\lambda,\mu})$. Then one also has that $\hat u_{\lambda,\mu}\in \arg\min_{u\in\realr^d} f(u)+\lambda (g(u)-\tau_1)+\mu(h(Au)-\tau_2)$ which implies
\begin{displaymath}
\begin{array}{lcl}
\tilde f(\lambda,\mu) & \stackrel{(\ref{eq:dual})}{=} & \inf_{u\in\realr^d}f(u)+\lambda (g(u)-\tau_1)+\mu(h(Au)-\tau_2) \\[3mm]
 & \stackrel{\text{(hypothesis)}}{=} & f(\hat u_{\lambda,\mu})+\lambda (g(\hat u_{\lambda,\mu})-\tau_1)+\mu(h(A\hat u_{\lambda,\mu})-\tau_2)\\[3mm]
& \stackrel{\text{(def. of }\tau )}{=} & f(\hat u_{\lambda,\mu})
\end{array}
\end{displaymath}
(strong duality for problems (\ref{eq:constrainedproblem}) and (\ref{eq:dual})) which implies that $\hat u_{\lambda,\mu}$ is optimal for (\ref{eq:constrainedproblem}) with $\tau=(g(\hat u_{\lambda,\mu}),h(A\hat u_{\lambda,\mu}))$ and $(\lambda,\mu)$ is optimal for the dual problem (\ref{eq:dual}).

5) As we have shown in point~\ref{item:duality} that $\hat u_{\lambda,\mu}$ is a minimizer of the constrained problem (\ref{eq:constrainedproblem}) for $\tau=(g(\hat u_{\lambda,\mu}),h(A\hat u_{\lambda,\mu}))$, it follows that the point with coordinates $(g(\hat u_{\lambda,\mu}),h(A\hat u_{\lambda,\mu}),f(\hat u_{\lambda,\mu}))$ lies on the Pareto frontier.

6) Let $\lambda,\mu\geq 0$. By point~\ref{item:duality} we know that if $\hat u_{\lambda,\mu}$ solves optimization problem $\min_{u\in\realr^d} f(u)+\lambda g(u)+\mu h(Au)$, then $\hat u_{\lambda,\mu}$ also solves the constrained problem (\ref{eq:constrainedproblem}) for $\tau=(g(\hat u_{\lambda,\mu}),h(A\hat u_{\lambda,\mu}))$ and strong duality holds. Thus one can write:
\begin{displaymath}
\begin{array}{lcl}
\varphi(\tau) &\stackrel{(\ref{eq:valuefunction})}{=}& \displaystyle \inf_{u\in S(\tau)} f(u)\\[3mm]
 &\stackrel{\text{strong duality}}{=}& \displaystyle \tilde f(\lambda,\mu)\\[3mm]
 &\stackrel{(\ref{eq:dual})}{\leq}& \displaystyle f(u)+\lambda (g(u)-\tau_1)+\mu (h(A u)-\tau_2)\qquad \forall u\in\realr^d\\[3mm]
\end{array}
\end{displaymath}
which implies $f(u)\geq\varphi(\tau)-\lambda (g(u)-\tau_1)-\mu(h(Au)-\tau_2)$. Therefore one must have:
\begin{displaymath}
\begin{array}{lcl}
\varphi(\tilde\tau)&\stackrel{\text{def.}}{=}&\displaystyle \inf_{u\in S(\tilde\tau)} f(u)\\[3mm]
& \geq & \varphi(\tau)+ \inf_{u\in S(\tilde\tau)}- \lambda(g(u)-\tau_1)-\mu(h(Au)-\tau_2)\\[3mm]
& \geq & \varphi(\tau)+ \inf_{u\in S(\tilde\tau)}- \lambda(g(u)-\tau_1)+ \inf_{u\in S(\tilde\tau)}-\mu(h(Au)-\tau_2)\\[3mm]
& \geq & \varphi(\tau)+ \inf_{g(u)\leq \tilde\tau_1}- \lambda(g(u)-\tau_1)+ \inf_{h(Au)\leq \tilde\tau_2}-\mu(h(Au)-\tau_2)\\[3mm]
& \geq & \varphi(\tau)-\lambda (\tilde\tau_1-\tau_1)-\mu (\tilde\tau_2-\tau_2)\qquad \forall \tilde \tau.
\end{array}
\end{displaymath}
This means that $-(\lambda,\mu)$ is an element of the subgradient of $\varphi$ at $\tau$: $-(\lambda,\mu)\in \partial \varphi(\tau)$.

7) Let us start by remarking that for $\tau$ in the relative interior of the domain of $\varphi$, one has $\partial \varphi (\tau)\neq \emptyset$ \cite[Th.~23.4]{Rockafellar1970}. In this case, set $-(\lambda,\mu)\in\partial \varphi(\tau)$. As $\varphi$ is decreasing we know that $\lambda,\mu\geq 0$. We can write $\varphi(\tilde\tau)\geq \varphi(\tau)- \lambda(\tilde\tau_1-\tau_1)-\mu(\tilde\tau_2-\tau_2)$ for all $\tilde\tau$.

Now suppose that $\tilde u_{\tau}$ is not a minimizer of $f(u)+\lambda g(u)+\mu h(Au)$, then there is a $\hat u\in\realr^d$ such that
\begin{displaymath}
f(\hat u)+\lambda g(\hat u)+\mu h(A\hat u) <  f(\tilde u_{\tau})+\lambda g(\tilde u_{\tau})+\mu h(A \tilde u_\tau)
\end{displaymath}
and thus
\begin{align}
f(\hat u) &<  f(\tilde u_{\tau})+\lambda \left[ g(\tilde u_{\tau})- g(\hat u)\right]
+\mu \left[ h(A\tilde u_{\tau})- h(A\hat u)\right]\nonumber\\[3mm]
& = \varphi(\tau)+\lambda \left[ g(\tilde u_{\tau})- g(\hat u)\right]
+\mu \left[ h(A\tilde u_{\tau})- h(A\hat u)\right]\nonumber\\[3mm]
& \leq \varphi(\tilde\tau)+\lambda (\tilde\tau_1-\tau_1)+\lambda \left[ g(\tilde u_{\tau})- g(\hat u)\right]
+\mu (\tilde\tau_2-\tau_2)+\mu \left[ h(A\tilde u_{\tau})- h(A\hat u)\right]\nonumber
\end{align}
($\forall \tilde\tau$). Now we choose $\tilde \tau=(g(\hat u),h(A \hat u))$ in the previous expression to find:
\begin{displaymath}
f(\hat u) <  \varphi(\tilde\tau)+\lambda (g(\tilde u_{\tau})-\tau_1)+\mu (h(A\tilde u_{\tau})-\tau_2) \leq \varphi(\tilde\tau)
\end{displaymath}
as $g(\tilde u_{\tau})\leq \tau_1$ and $h(A\tilde u_{\tau})\leq \tau_2$ by definition. We have by definition of $\tilde \tau$ that $\hat u$ satisfies $(g(\hat u),h(A\hat u))\leq \tilde \tau$ and we have shown that $f(\hat u)<\varphi(\tilde\tau)$ which is in contradiction to the definition of $\varphi$. Therefore $\tilde u_{\tau}$ must indeed be a minimizer of $f(u)+\lambda g(u)+\mu h(Au)$.
\end{proof}

Taken together, the seven points of Theorem~\ref{prop:Pareto} justify the description of the Pareto curve for problem (\ref{eq:problem}) given in the Introduction and clarify the role of $-(\lambda,\mu)$ as subgradient to it.

\section{A primal-dual continuation algorithm}\label{sec-pd-alg}

In this section we modify the primal-dual algorithm of \cite{Condat2013} to include variable penalty parameters $\lambda_n$ and $\mu_n$, and show its convergence to a minimizer of the cost function in problem (\ref{eq:problem}). As mentioned in the introduction, if the sequences $(\lambda_n)_{n\in\naturaln}$ and $(\mu_n)_{n\in\naturaln}$ are chosen to vary slowly with $n$, one may expect that all the iterates $(u_n)_{n\in\naturaln}$ provide points $((g(u_n),h(Au_n)),f(u_n))$ near the Pareto frontier of the problem. Such a strategy is therefore advantageous if one wants to construct (a part of) this frontier. 

The proof of convergence of algorithm (\ref{eq:algorithm}) given here does not rely on the reformulation in terms of monotone operators technique used in \cite{Condat2013}, but it is more explicit, as it will involve the use of a variable metric to measure the convergence to a fixed point.

We start by writing the variational equations that the minimizers $\hat u$ of the convex optimization problem (\ref{eq:problem}) must satisfy:
\begin{displaymath}
     \grad f(\hat u)+ \hat w+\mu A^T \hat v=0 \qquad\text{with}\qquad
    \hat w\in\partial\lambda g(\hat u)\quad \text{and}\quad \hat v\in\partial h(A\hat u).
\end{displaymath}
Using the relations (\ref{eq:subdiffcharac}), this can also be written as
\begin{displaymath}
\left\{
\begin{array}{l}
\grad f(\hat u)+\hat w+\mu A^T \hat v=0,\\% \quad
\hat u=\prox_{\alpha \lambda g}(\hat u+\alpha \hat w),\\% \quad\text{and}\quad
\hat v=\prox_{(\beta/\mu) h^\ast}(\hat v+(\beta/\mu) A\hat u),
\end{array}
\right.
\end{displaymath}
for any $\alpha,\beta>0$ or by eliminating the auxiliary variable $\hat w$:
\begin{align}
\hat u &=\prox_{\alpha\lambda g}(\hat u-\alpha \grad f(\hat u)-\alpha\mu A^T \hat v) \label{eq:fixedpointa}\\
\hat v &=\prox_{(\beta/\mu) h^\ast}(\hat v+(\beta/\mu) A\hat u)\label{eq:fixedpointb}
\end{align}
which are now conveniently written as fixed-point equations.
In order to solve the equations (\ref{eq:fixedpointa})--(\ref{eq:fixedpointb}) and hence the optimization problem (\ref{eq:problem}), one can use the iterative algorithm proposed by \cite{Condat2013} which takes the form:
\begin{equation}\label{eq:alg1}
\left\{
\begin{array}{lcl}
u_{n+1} &=& \prox_{\alpha\lambda g}\left(u_n-\alpha \grad f(u_n)-\alpha\mu A^T  v_n\right)\\
v_{n+1} &=& \prox_{(\beta/\mu)\, h^\ast}\left(v_n+(\beta/\mu) A (2u_{n+1}-u_n)\right)
\end{array}
\right.
\end{equation}
for any $u_0\in \realr^d,v_0\in\realr^{d'}$. In \cite{Condat2013} it is shown that this algorithm converges as long as the step length parameters are chosen to satisfy $\beta \|A\|^2<1/\alpha-L/2$ where $L$ is a Lipschitz constant of the gradient of $f$. 
%In fact it is also possible to formulate this algorithm with additional error terms to take into account numerical errors in computing the right hand sides \cite{Condat2013}.

As announced, in this paper we propose a modification of algorithm (\ref{eq:alg1}) which allows for the use of variable penalty parameters $\lambda_n $ and $\mu_n$ instead of $\lambda,\mu$ fixed (independent of $n$): 
\begin{align}
& u_{n+1} = \prox_{\alpha\lambda_n g}\left(u_n-\alpha \grad f(u_n)-\alpha\mu_n A^T  v_n\right)\label{eq:alg1a}\\
& v_{n+1} = \prox_{(\beta/\mu_n)\, h^\ast}\left(v_n+(\beta/\mu_n) A (2u_{n+1}-u_n)\right)\label{eq:alg1b}
\end{align}
with $(\lambda_n)_{n\in\naturaln}, (\mu_n)_{n\in\naturaln}\subset \realr_{>0}$ two a priori sequences converging to the penalty parameters $\lambda$ and $\mu$ of problem (\ref{eq:problem}):  $\lambda_n\stackrel{n\to\infty}{\longrightarrow}\lambda$ and $\mu_n\stackrel{n\to\infty}{\longrightarrow}\mu$. The proof of convergence relies on the following lemmas.

\begin{lemma} If $f:\realr^d \rightarrow \realr$ is convex with Lipschitz continuous gradient (constant $L$) then $L^{-1}\grad f$ is firmly non expansive:
\begin{equation}\label{fne}
\langle L^{-1}\grad f(u)-L^{-1}\grad f(v),u-v\rangle \geq \|L^{-1}\grad  f(u)-L^{-1}\grad f(v)\|_2^2
%\quad\forall 
\end{equation}
for all $u,v\in\realr^d$.
\end{lemma}
\begin{proof}
See \cite[Part 2, Chapter X, Th. 4.2.2]{HiriartUrruty1993}.
%\stopproof
\end{proof}

\begin{lemma}
    Let $g\in \Gamma_0(\realr^d)$, then $u_+=\prox_{g}(u_-+\Delta)$ if and only if
    \begin{equation}\label{eq:prox_ineq}
        \|u_+-u\|_2^2\leq \|u_--u\|_2^2-\|u_+-u_-\|_2^2+2\langle u_+-u,\Delta\rangle+ 2g(u)-2g(u_+) 
    \end{equation}
for all $\forall \ u\in\realr^d$.
\end{lemma}
\begin{proof}
The relation $u_+=\prox_{g}(u_-+\Delta)$ means that
$u_+=\arg\min_u \frac{1}{2}\|u-(u_-+\Delta)\|_2^2 + g(u)$ or
$0\in u_+-(u_-+\Delta)+\partial g(u_+)$. The inclusion $u_--u_++\Delta\in \partial g(u_+)$ is equivalent to:
\begin{displaymath}
g(u)\geq g(u_+)+\langle u_--u_++\Delta,u-u_+\rangle\qquad\forall u\in\realr^d.
\end{displaymath}
The inner product can be re-arranged to give relation (\ref{eq:prox_ineq}).
%\stopproof
\end{proof}

\begin{lemma}\label{lemma-inequality}
Let $\{a_n\}_{n\in\naturaln},\{b_n\}_{n\in\naturaln}\subset \realr_{\geq 0}$ be two summable sequences:  $\sum_n a_n<\infty$ and $\sum_n b_n<\infty$. If the sequence $\{\epsilon_n\}_{n\in\naturaln}\subset \realr_{\geq 0}$ satisfies $\epsilon_{n+1}^2\leq (1+2a_n)\epsilon_n^2+2 b_n \epsilon_{n+1}$ then it is bounded.
\end{lemma}    
\begin{proof} It is a slight generalization of \cite[Section 2.2.1, Lemma 2]{Polyak1987}.
We start by rewriting the inequality  as: $(\epsilon_{n+1}-b_n)^2\leq (1+2a_n) \epsilon_n^2+b_n^2$ which therefore implies 
\begin{displaymath}
|\epsilon_{n+1}-b_n|\leq \sqrt{(1+2a_n) \epsilon_n^2+b_n^2}\leq \sqrt{1+2a_n}\epsilon_n+b_n\leq (1+a_n)\epsilon_n+b_n.
\end{displaymath}
In case $\epsilon_{n+1}-b_n\geq 0$ one finds $\epsilon_{n+1}\leq (1+a_n)\epsilon_n+2b_n$.
In case $\epsilon_{n+1}-b_n\leq 0$ one also finds $\epsilon_{n+1}\leq b_n\leq (1+a_n)\epsilon_n+2b_n$.
Finally, this implies 
\begin{displaymath}
\begin{array}{lcl}
\epsilon_{n+1} & \leq  & \displaystyle (1+a_n)\epsilon_n+2b_n \\[4mm]
 & \leq & \displaystyle \prod_{k=0}^n(1+a_k)\epsilon_0+2\sum_{k=0}^n\prod_{l=k+1}^n(1+a_l) b_k\\[5mm]
 & \leq & \displaystyle \prod_{k=0}^\infty(1+a_k)\epsilon_0+2\sum_{k=0}^\infty\prod_{l=0}^\infty(1+a_l)b_k\\[5mm]
&<& \infty
\end{array}
\end{displaymath}
independently of $n$.
\end{proof}

\begin{lemma}\label{lemma:norm}
Let $\mu>0, L\geq 0$, $\{\mu_n\}_{n\in\naturaln}\subset \realr_{>0}$ such that $\sum_n |\mu_n-\mu|<\infty$. Let $A:\realr^d\to\realr^{d'}$ be a linear map and $\alpha,\beta>0$ such that $\beta\|A\|^2<\alpha^{-1}-L/2$, then the matrices
\begin{displaymath}
M_n=\left(
\begin{array}{cc}
\mu_n^{-1}\alpha^{-1}\id & -A^T \\
-A & \mu_n\beta^{-1}\id
\end{array}
\right)\quad\text{and}\quad 
M=\left(
\begin{array}{cc}
\mu^{-1}\alpha^{-1}\id & -A^T \\
-A & \mu\beta^{-1}\id
\end{array}
\right)
\end{displaymath}
are positive definite and define norms 
\begin{displaymath}
\|(u,v)\|_{M_n}=\sqrt{\langle (u,v), M_n(u,v)\rangle}\quad\text{and}\quad \|(u,v)\|_{M}=\sqrt{\langle (u,v), M(u,v)\rangle}
\end{displaymath}
on $\realr^{d+d'}$. Furthermore, one has the bounds
\begin{equation}\label{eq:normbound}
\begin{array}{l}
\|(u,v)\|_{M_n}^2\geq\left\|(u,v)\right\|_{M_n}^2-L\|u\|_2^2/(2\mu_n)\geq C_1 \|(u,v)\|_2^2\\[3mm]%\quad\text{and}\quad
 \|(u,v)\|_{M}\geq C_2 \|(u,v)\|_2
\end{array}
\end{equation}
for some $C_1, C_2>0$ independent of $n$.
There also exist summable sequences $(e_n)_{n\in\naturaln},(f_n)_{n\in\naturaln}\subset\realr_{>0}$ such that 
\begin{equation}\label{eq:norm-equiv}
\|x\|_{M}\leq (1+e_n)\|x\|_{M_n}\quad\text{and}\quad \|x\|_{M_n}\leq (1+f_n)\|x\|_{M}
\end{equation}
for all $x\in\realr^{d+d'}$.
\end{lemma}
\begin{proof} We first write $2\langle Au,v\rangle=2\langle \sqrt{\epsilon/\mu_n}\, u,\sqrt{\mu_n/\epsilon}A^T v\rangle\leq (\epsilon/\mu_n) \|u\|_2^2+(\mu_n/\epsilon) \|A^T v\|_2^2$ for any $\epsilon>0$. Therefore we have 
\begin{displaymath}
\begin{array}{lcl}
\langle (u,v), M_n(u,v)\rangle -L\|u\|_2^2/2\mu_n& = & \mu_n^{-1}\alpha^{-1}\|u\|_2^2+\mu_n\beta^{-1}\|v\|_2^2 -2 \langle Au,v\rangle\\[4mm]
&&\qquad\qquad\qquad\qquad\qquad\displaystyle  -L\|u\|_2^2/2\mu_n\\[4mm]
&\geq &  \mu_n^{-1}(\alpha^{-1}-L/2-\epsilon)\|u\|_2^2\\[4mm]
&&\displaystyle\qquad \qquad+\mu_n(\beta^{-1}-\epsilon^{-1} \|A\|^2) \|v\|_2^2.
\end{array}
\end{displaymath}
We have by virtue of the hypothesis $\beta \|A\|^2<\alpha^{-1}-L/2$, that there exists $\epsilon>0$ and $\delta>0$ (small) such that $[\epsilon-\delta\beta\epsilon,\epsilon+\delta]\subset (\beta \|A\|^2, \alpha^{-1}-L/2)$ which implies
\begin{displaymath}
\alpha^{-1}-L/2-\epsilon >\delta \quad\text{and}\quad \beta^{-1}-\epsilon^{-1}\|A\|^2>\delta.
\end{displaymath}
Then we continue from above:
\begin{displaymath}
\begin{array}{lcl}
\langle (u,v), M_n(u,v)\rangle-L\|u\|_2^2/2\mu_n &\geq &  \mu_n^{-1}(\alpha^{-1}-L/2-\epsilon)\|u\|_2^2\\[4mm]
&& \displaystyle \qquad \qquad+\mu_n(\beta^{-1}-\epsilon^{-1} \|A\|^2) \|v\|_2^2\\[4mm]
&\geq  & \mu_n^{-1}\delta \|u\|_2^2 +\mu_n\delta \|v\|_2^2\\[4mm]
&\geq& C_1 \|(u,v)\|_2^2
\end{array}
\end{displaymath}
where $C_1>0$ does not depend on $n$. On the one hand this shows that $M_n$ is positive definite and on the other hand it shows that inequalities (\ref{eq:normbound}) hold. One proves in the same way that $M$ is also positive definite.

In order to show inequalities (\ref{eq:norm-equiv}) we first write 
\begin{displaymath}
M_n=M+
(\mu_n-\mu)\Delta_n\quad\text{with}\quad \Delta_n=\left(
\begin{array}{cc}
-(\alpha\mu\mu_n)^{-1}\id & 0\\
0 & \beta^{-1}\id
\end{array}
\right)
\end{displaymath}
with $\|\Delta_n\|=\max((\alpha\mu\mu_n)^{-1},\beta^{-1})\leq C_3$ (independent of $n$) and hence find:
\begin{displaymath}
\begin{array}{lcl}
\|x\|_{M}^2&=&\|x\|_{M_n}^2-(\mu_n-\mu)\langle x, \Delta_n x\rangle\\[4mm]
&\leq &  \|x\|_{M_n}^2+|\mu_n-\mu|\, \|\Delta_n\|\, \|x\|_2^2\\[4mm]
&\stackrel{(\ref{eq:normbound})}{\leq}& (1+C_3 C_1^{-1}|\mu_n-\mu|)\|x\|_{M_n}^2. 
\end{array}
\end{displaymath}
such that $\|x\|_{M}\leq \sqrt{1+C_3 C_1^{-1}|\mu_n-\mu|}\times \|x\|_{M_n}\leq (1+C_3 C_1^{-1}|\mu_n-\mu|/2)\|x\|_{M_n}$.
On the other hand, one also has
\begin{displaymath}
\begin{array}{lcl}
\|x\|_{M_n}^2&=&\|x\|_{M}^2+(\mu_n-\mu)\langle x, \Delta_n x\rangle\\[4mm]
&\leq&  \|x\|_{M}^2+|\mu_n-\mu|\, \|\Delta_n\|\, \|x\|_2^2\\[4mm]
&\stackrel{(\ref{eq:normbound})}{\leq}& (1+C_3/C_2 |\mu_n-\mu|)\|x\|_{M}^2. 
\end{array}
\end{displaymath}
Therefore we find $\|x\|_{M_n}\leq \sqrt{1+C_3/C_2 |\mu_n-\mu|}\times \|x\|_{M} \leq (1+C_3/C_2 |\mu_n-\mu|/2)\|x\|_{M}$ as well.
\end{proof}

\begin{theorem}\label{prop-convergence}
If a minimizer $\hat u$ of problem (\ref{eq:problem}) (with the stated assumptions) exists and the sequence $(u_n)_{n\in\naturaln}$ is defined  by the algorithm (\ref{eq:alg1a})--(\ref{eq:alg1b}) with $\alpha,\beta>0$ and $\beta \|A\|^2<1/\alpha-L/2$ and with parameters $(\lambda_n)_{n\in\naturaln}, (\mu_n)_{n\in\naturaln}\subset\realr_{>0}$ satisfying
\begin{equation}
\sum_{n\in\naturaln}|\lambda_n-\lambda|<\infty \quad\text{and}\quad \sum_{n\in\naturaln}|\mu_n-\mu|<\infty 
\end{equation}
then the sequence $(u_n)_{n\in\naturaln}$ converges to a minimizer of problem (\ref{eq:problem}) for any starting point $(u_0,v_0)$.
\end{theorem}
\begin{proof} Let $\hat u\in \arg\min_u f(u)+\lambda g(u)+\mu h(Au)$, i.e. $\hat u$ and $\hat v$ satisfy equations (\ref{eq:fixedpointa}) and (\ref{eq:fixedpointb}). We first apply inequality (\ref{eq:prox_ineq}) to the equation (\ref{eq:alg1a}) with $u^+=u_{n+1}, u^-=u_n, u=\hat u$ and $\Delta=-\alpha \grad f(u_n)-\alpha\mu_n A^Tv_n$:
\begin{displaymath}
\begin{array}{lcl}
\|u_{n+1}-\hat u\|_2^2&\leq& \|u_n-\hat u\|_2^2-\|u_{n+1}-u_n\|_2^2\\[2mm]
&& \qquad\qquad\qquad
+2\langle u_{n+1}-\hat u,-\alpha \grad f(u_n)-\alpha\mu_n A^T v_n\rangle\\[2mm]
&& \qquad\qquad\qquad\qquad\qquad\qquad+2\alpha\lambda_n g(\hat u)-2\alpha\lambda_n g(u_{n+1})
\end{array}
\end{displaymath}
and to the fixed point relation (\ref{eq:fixedpointa}) with  $u^+=\hat u, u^-=\hat u, u=u_{n+1}$  and $\Delta=-\alpha \grad f(\hat u)-\alpha\mu A^T\hat v$:
\begin{displaymath}
\begin{array}{lcl}
\|\hat u- u_{n+1}\|_2^2&\leq& \|\hat u - u_{n+1}\|_2^2-\|\hat u-\hat u\|_2^2\\[2mm]
&&\qquad\qquad\qquad+2\langle \hat u -u_{n+1},-\alpha \grad f(\hat u)-\alpha\mu  A^T \hat v\rangle\\[2mm]
&& \qquad\qquad\qquad\qquad\qquad\qquad+2\alpha\lambda g( u_{n+1})-2\alpha\lambda g(\hat u).
\end{array}
\end{displaymath}
Summing the first with $\lambda_n/\lambda$ times the second inequality gives:
\begin{equation}\label{eq:tmp1}
\begin{array}{lcl}
\|u_{n+1}-\hat u\|_2^2&\leq& \|u_n-\hat u\|_2^2-\|u_{n+1}-u_n\|_2^2 \\[2mm]
&&\qquad-2\alpha\langle u_{n+1}-\hat u, \grad f(u_n) -\grad f(\hat u)\rangle\\[2mm]
&&\qquad-2\alpha(\lambda-\lambda_n)/\lambda \langle u_{n+1}-\hat u,  \grad f(\hat u)\rangle\\[2mm]
 &&\qquad-2\alpha\mu_n \langle u_{n+1}-\hat u, A^T (v_n-\hat v)\rangle \\[2mm]
&&\qquad-2\alpha/\lambda \langle u_{n+1}-\hat u, (\mu_n\lambda-\mu\lambda_n)A^T \hat v\rangle.
\end{array}
\end{equation}
The inner product $\langle\hat u- u_{n+1},\grad f(u_n)-\grad f(\hat u)\rangle$ can be bounded by $L\|u_{n+1}-u_n\|_2^2/4$ because:
\begin{displaymath}
\begin{array}{l}
\displaystyle-\langle u_{n+1}-\hat u ,\grad f(u_n)-\grad f(\hat u)\rangle\\[3mm]
\qquad = \displaystyle\langle\hat u- u_{n},\grad f(u_n)-\grad f(\hat u)\rangle+\langle u_n- u_{n+1},\grad f(u_n)-\grad f(\hat u)\rangle\\[3mm]
\qquad \stackrel{(\ref{fne})}{\leq}\displaystyle\frac{-1}{L}\|\grad f(u_n)-\grad f(\hat u)\|_2^2+\langle u_n- u_{n+1},\grad f(u_n)-\grad f(\hat u)\rangle\\[3mm]
\qquad = \displaystyle\langle \sqrt{L}(u_n- u_{n+1})-\frac{1}{\sqrt{L}}(\grad f(u_n)-\grad f(\hat u)),\frac{1}{\sqrt{L}}(\grad f(u_n)-\grad f(\hat u))\rangle\\[3mm]
\qquad = \displaystyle \frac{L}{4}\left\|u_n- u_{n+1}\right\|_2^2-\frac{1}{4}\left\|u_n- u_{n+1}-\frac{2}{\sqrt{L}}(\grad f(u_n)-\grad f(\hat u))\right\|_2^2\\[5mm]
\qquad \leq \displaystyle \frac{L}{4}\|u_n- u_{n+1}\|_2^2.
\end{array}
\end{displaymath}
Therefore, using also $\mu_n\lambda-\mu\lambda_n=\mu(\lambda-\lambda_n)+\lambda(\mu_n-\mu)$, inequality (\ref{eq:tmp1}) implies:
\begin{equation}\label{eq:tmp2}
\begin{array}{lcl}
\|u_{n+1}-\hat u\|_2^2 & \leq & \|u_n-\hat u\|_2^2-(1-\alpha L/2)\|u_{n+1}-u_n\|_2^2\\[2mm]
&& \qquad\qquad\quad-2\alpha\mu_n \langle u_{n+1}-\hat u, A^T (v_n-\hat v)\rangle\\[2mm]
 &&\qquad\qquad\quad + (c_1 |\lambda_n-\lambda|+c_2|\mu_n-\mu|)\|u_{n+1}-\hat u\|_2
\end{array}
\end{equation}
with $c_1,c_2\in\realr_{>0}$ independent of $n$.

Now we apply inequality (\ref{eq:prox_ineq}) to the relation (\ref{eq:alg1b}) with $u^+=v_{n+1}, u^-=v_n, u=\hat v$ and $\Delta=(\beta/\mu_n)A(2u_{n+1}-u_n)$ to give
\begin{displaymath}
\begin{array}{lcl}
\|v_{n+1}-\hat v\|_2^2&\leq& \|v_n-\hat v\|_2^2-\|v_{n+1}-v_n\|_2^2 +2(\beta/\mu_n)\langle v_{n+1}-\hat v,  A(2u_{n+1}-u_n)\rangle
\\[2mm]
&&\qquad + 2(\beta/\mu_n) h^\ast(\hat v)-2(\beta/\mu_n) h^\ast(v_{n+1})
\end{array}
\end{displaymath}
and likewise we apply  inequality (\ref{eq:prox_ineq}) to the relation (\ref{eq:fixedpointb}) with $u^+=\hat v, u^-=\hat v, u=v_{n+1}$ and $\Delta=\beta/\mu A\hat u$:
\begin{displaymath}
\begin{array}{lcl}
\|\hat v-v_{n+1}\|_2^2&\leq& \|\hat v-v_{n+1}\|_2^2-\|\hat v-\hat v\|_2^2+2(\beta/\mu)\langle \hat v-v_{n+1}, A\hat u\rangle\\[2mm]
&&\qquad +2(\beta/\mu) h^\ast(v_{n+1})-2(\beta/\mu) h^\ast(\hat v).
\end{array}
\end{displaymath}
Now we add $\mu/\mu_n$ times the latter to the former so as to obtain:
\begin{equation}\label{eq:tmp3}
\|v_{n+1}-\hat v\|_2^2\leq \|v_n-\hat v\|_2^2-\|v_{n+1}-v_n\|_2^2+2(\beta/\mu_n)\langle v_{n+1}-\hat v, A(2u_{n+1}-u_n-\hat u)\rangle
\end{equation}
Summing $1/(\alpha\mu_n)$ times inequality (\ref{eq:tmp2}) and $(\mu_n/\beta)$ times inequality (\ref{eq:tmp3}) gives:
\begin{displaymath}%\label{eq:tmp4}
\begin{array}{l}
\|u_{n+1}-\hat u\|_2^2/(\alpha\mu_n)+ (\mu_n/\beta)\|v_{n+1}-\hat v\|_2^2 \leq \\[3mm]
\qquad\qquad\qquad\qquad\qquad\qquad \|u_n-\hat u\|_2^2/\alpha\mu_n+(\mu_n/\beta)\|v_n-\hat v\|_2^2\\[3mm]
\qquad\qquad\qquad\qquad\qquad\qquad -(1/\alpha-L/2)/\mu_n \|u_{n+1}-u_n\|_2^2-(\mu_n/\beta)\|v_{n+1}-v_n\|_2^2\\[3mm]
\qquad\qquad\qquad\qquad\qquad\qquad -2\langle u_{n+1}-\hat u, A^T (v_n-\hat v)\rangle\\[3mm]
\qquad\qquad\qquad\qquad\qquad\qquad +2\langle v_{n+1}-\hat v, A(2u_{n+1}-\hat u-u_n)\rangle\\[3mm]
\qquad\qquad\qquad\qquad\qquad\qquad + (c_1 |\lambda_n-\lambda|+c_2|\mu_n-\mu|)\|u_{n+1}-\hat u\|_2
\end{array}
\end{displaymath}
for some (new) constants $c_1,c_2>0$.
The two remaining scalar products can be written as: 
\begin{displaymath}
\begin{array}{lcl}
2\langle u_{n+1}-\hat u,- A^T (v_n-\hat v)\rangle+2\langle v_{n+1}-\hat v, A(2u_{n+1}-\hat u-u_n)\rangle\\[2mm]
\qquad\qquad\qquad\qquad=2\langle u_{n+1}-\hat u, A^T (v_{n+1}-\hat v)\rangle-2\langle u_n-\hat u,A^T (v_n-\hat v)\rangle\\[2mm]
\qquad\qquad\qquad\qquad\qquad\qquad+2\langle u_{n+1}-u_n,A^T (v_{n+1}-v_n)\rangle,
\end{array}
\end{displaymath}
hence the previous inequality reduces to:
\begin{displaymath}
\begin{array}{l}
\|u_{n+1}-\hat u\|_2^2/\alpha\mu_n+(\mu_n/\beta) \|v_{n+1}-\hat v\|_2^2 -2\langle u_{n+1}-\hat u, A^T (v_{n+1}-\hat v)\rangle\\[2mm]
\qquad \leq \|u_n-\hat u\|_2^2/\alpha\mu_n+\mu_n\|v_n-\hat v\|_2^2/\beta -2\langle u_n-\hat u,A^T (v_n-\hat v)\rangle\\[2mm]
\qquad -(1/\alpha-L/2)/\mu_n\|u_{n+1}-u_n\|_2^2-(\mu_n/\beta)\|v_{n+1}-v_n\|_2^2\\[2mm]
\qquad +2\langle u_{n+1}-u_n,A^T (v_{n+1}-v_n)\rangle + (c_1 |\lambda_n-\lambda|+c_2|\mu_n-\mu|)\|u_{n+1}-\hat u\|_2,
\end{array}
\end{displaymath}
which can be rewritten in terms of the norms $\|\,\,\,\|_{M_n}$ introduced in Lemma~\ref{lemma:norm} as
\begin{displaymath}
\begin{array}{lcl}
\left\|\stackvec{u_{n+1}-\hat u}{v_{n+1}-\hat v}\right\|_{M_n}^2
&\leq& \left\|\stackvec{u_n-\hat u}{v_n-\hat v}\right\|_{M_n}^2
-\left\|\stackvec{u_{n+1}-u_n}{v_{n+1}-v_n}\right\|_{M_n}^2\\[6mm]
&&\qquad +L\|u_{n+1}-u_n\|_2^2/2\mu_n\\[4mm]
&&\qquad + (c_1 |\lambda_n-\lambda|+c_2|\mu_n-\mu|)\|u_{n+1}-\hat u\|_2.
\end{array}
\end{displaymath}
Setting $\hat x=(\hat u,\hat v)$, $x_{n+1}=(u_{n+1},v_{n+1})$, $x_n=(u_n,v_n)$ and using inequality (\ref{eq:normbound}) one finds:
\begin{equation}\label{eq:tmp5}
\begin{array}{lcl}
\|x_{n+1}-\hat x\|_{M_n}^2 &\leq& \|x_{n}-\hat x\|_{M_n}^2 -C_1 \left\|x_{n+1}-x_{n}\right\|_2^2 \\[2mm]
 &&\qquad\qquad+ (\bar c_1 |\lambda_n-\lambda|+\bar c_2|\mu_n-\mu|)\left\|x_{n+1}-\hat x\right\|_{M}
\end{array}
\end{equation}
where we have also used the bounds $\|u_{n+1}-\hat u\|_2\leq \|x_{n+1}-\hat x\|_2\leq  \|x_{n+1}-\hat x\|_M/C_2$ and have introduced new constants $\bar{c}_1=c_1/C_2$ and $\bar{c}_2=c_2/C_2$.

One finds:
\begin{displaymath}
\begin{array}{lcl}
\|x_{n+1}-\hat x\|_{M}^2 & \stackrel{(\ref{eq:norm-equiv})}{\leq} & (1+e_n) \|x_{n+1}-\hat x\|_{M_n}^2\\[3mm]
& \stackrel{(\ref{eq:tmp5})}{\leq} & (1+e_n) \left[\rule{0pt}{5mm}  \|x_{n}-\hat x\|_{M_n}^2 -C_1 \left\|x_{n+1}-x_{n}\right\|_2^2 \right. \\[3mm]
&& \qquad\qquad \left. \rule{0pt}{5mm}+ (\bar c_1 |\lambda_n-\lambda|+\bar c_2|\mu_n-\mu|)\left\|x_{n+1}-\hat x\right\|_{M}\right] \\[3mm]
& \stackrel{(\ref{eq:norm-equiv})}{\leq} & (1+e_n)(1+f_n) \|x_{n}-\hat x\|_{M}^2 -\tilde C_1 \left\|x_{n+1}-x_{n}\right\|_2^2  \\[3mm]
&& \qquad\qquad  + (\tilde c_1 |\lambda_n-\lambda|+\tilde c_2|\mu_n-\mu|)\left\|x_{n+1}-\hat x\right\|_{M} \\[3mm]
\end{array}
\end{displaymath}
or
\begin{equation}\label{eq:fund-ineq}
\|x_{n+1}-\hat x\|_{M}^2  \leq  (1+2 a_n) \|x_{n}-\hat x\|_{M}^2 -c \left\|x_{n+1}-x_{n}\right\|_2^2    + 2b_n \left\|x_{n+1}-\hat x\right\|_{M} \\[3mm]
\end{equation}
with positive $a_n,b_n, c$ such that $\sum_na_n<\infty$ and $\sum_n b_n<\infty$.

It follows from applying Lemma~\ref{lemma-inequality} to inequality (\ref{eq:fund-ineq}) that the sequence $(u_n,v_n)_{n\in\naturaln}\subset\realr^{d+d'}$ is bounded. This implies the existence of an accumulation point $(u^\dagger,v^\dagger)$ and converging subsequences $u_{n_j}\stackrel{j\to \infty}{\longrightarrow} u^\dagger$ and $v_{n_j}\stackrel{j\to \infty}{\longrightarrow} v^\dagger$.

By summing the inequalities (\ref{eq:fund-ineq}) from $n=0$ to $n=N$, it also follows that 
\begin{displaymath}
\begin{array}{lcl}
\displaystyle c \sum_{n=0}^N \left\|x_{n+1}-x_{n}\right\|_2^2 & \leq & \displaystyle \sum_{n=0}^N  -\|x_{n+1}-\hat x\|_{M}^2+   (1+2 a_n) \|x_{n}-\hat x\|_{M}^2  \\[3mm]
&& \qquad\qquad\qquad\qquad\qquad\qquad\qquad    + 2b_n \left\|x_{n+1}-\hat x\right\|_{M}\\[3mm]
 & = & \displaystyle \|x_{0}-\hat x\|_{M}^2-\|x_{N+1}-\hat x\|_{M}^2+ \sum_{n=0}^N  2 a_n \|x_{n}-\hat x\|_{M}^2    \\[3mm]
&&\qquad\qquad\qquad\qquad\qquad\qquad\qquad + 2b_n \left\|x_{n+1}-\hat x\right\|_{M}\\[3mm]
&\leq & C 
\end{array}
\end{displaymath}
independently of $N$, as $(x_n)_{n\in\naturaln}$ is bounded. And thus $\lim_{n\to\infty}\|x_{n+1}-x_n\|_2=0$. It then follows from the continuity of expressions (\ref{eq:alg1a}) and (\ref{eq:alg1b}) that $(u^\dagger,v^\dagger)$ is a solution of equations (\ref{eq:fixedpointa}) and (\ref{eq:fixedpointb}). Hence we may replace $\hat x=(\hat u,\hat v)$ by $x^\dagger=(u^\dagger,v^\dagger)$ in relation (\ref{eq:fund-ineq}) and sum from $n=N_0$ to $n=N$ to obtain:
\begin{equation}
\begin{array}{lcl}
\|x_{N+1}- x^\dagger\|_{M}^2 & \leq & \displaystyle \|x_{N_0}-x^\dagger\|_{M}^2+ \sum_{n=N_0}^N  2 a_n \|x_{n}-x^\dagger\|_{M}^2     + 2b_n \left\|x_{n+1}-x^\dagger\right\|_{M}\\[3mm]
& \leq & \displaystyle \|x_{N_0}-x^\dagger\|_{M}^2+ \sum_{n=N_0}^N  \tilde a_n\\[3mm]
\end{array}
\end{equation}
(as $(x_n)_{n\in\naturaln}$ is bounded) where the sequence $(\tilde a_n)_{n_\in\naturaln}$ is again summable. Finally this implies that $\|x_{N+1}-x^\dagger\|_{M}<\epsilon$ for $N>N_0$ and $N_0$ sufficiently large on account of the converging subsequence and the summability of $(\tilde a_n)_{n\in\naturaln}$.
\end{proof}

There exist several special cases.
\begin{itemize}
\item The case $A=0$: The algorithm reduces to the proximal gradient algorithm 
\begin{displaymath}
u_{n+1}=\prox_{\alpha\lambda_n g}(u_n-\alpha \grad f(u_n))
\end{displaymath}
for the problem $\arg\min_u f(u)+\lambda g(u)$. A proof of convergence of this proximal-gradient-like continuation method was presented in \cite{Fest2023}.

\item In case $f=0$ one finds the algorithm
\begin{equation}
\left\{
\begin{array}{l}
u_{n+1} = \prox_{\alpha \lambda_n g}\left(u_n-\alpha \mu_n A^T  v_n\right)\\
v_{n+1} = \prox_{(\beta/\mu_n)\, h^\ast}\left(v_n+(\beta/\mu_n) A (2u_{n+1}-u_n)\right)
\end{array}
\right.
\end{equation}
which converges, when $\alpha\beta \|A\|^2<1$, to a minimizer of the problem:
\begin{displaymath}
\arg\min_u \lambda g(u)+\mu h(Au).
\end{displaymath}
It is the algorithm of Chambolle-Pock \cite{Chambolle2011}, but with variable parameters $\lambda_n$ and $\mu_n$.
\end{itemize}
In the same manner, one could also prove the convergence of the algorithm:
\begin{displaymath}
\left\{
\begin{array}{lcl}
v_{n+1} = \prox_{(\beta/\mu_n)\, h^\ast}\left(v_n+(\beta/\mu_n) A u_{n}\right)\\
u_{n+1} = \prox_{\alpha\lambda_n g}\left(u_n-\alpha \grad f(u_n)-\alpha\mu_n A^T  (2v_{n+1}-v_n)\right)
\end{array}
\right.
\end{displaymath}
under the same conditions as in Theorem~\ref{prop-convergence}.

\section{Relation to inexact proximal evaluations}\label{sec:inexact}

Under some additional assumptions on the functions $g$ and $h$, we can interpret method (2) as a special instance of the primal--dual method with inexact proximal evaluations proposed in \cite{Rasch2020}. Indeed, let us assume that $g$ is continuous on its domain and $h^*$ is continuous on a bounded and closed domain; the latter assumption is satisfied if, for instance, $h$ is any vectorial norm on $\mathbb{R}^{d'}$ \cite[Chapter 4]{Beck2017}. 

Let $\epsilon>0$. We recall that the $\epsilon-$subdifferential of a proper, convex, and closed function $F:\realr^d\rightarrow \realrplusinfty$ at the point $u\in\operatorname{dom}(F)$ is defined as \cite[Definition 1.1.1]{HiriartUrruty1993b}
\begin{equation}\label{eq:eps_sub}
\partial_{\epsilon}F(u)=\{\xi\in\realr^d: \ F(w)\geq F(u)+\langle \xi,w-u\rangle-\epsilon, \ \forall \ w\in\operatorname{dom}(F)\}.
\end{equation}
For any $a\in\mathbb{R}^d$, we introduce the function
\begin{equation*}
\varphi^{(a)}(u)=\frac{1}{2}\|u-a\|_2^2+F(u), \quad \forall \ u\in\mathbb{R}^d.
\end{equation*} 
Let $\hat{u}\in\realr^d$ be the (exact) proximal point of $F $ evaluated at point $a$, i.e., 
\begin{equation*}
\hat{u}=\operatorname{prox}_{F}(a) \quad \Leftrightarrow \quad \hat{u}=\underset{u\in\realr^d}{\operatorname{argmin}} \ \varphi^{(a)}(u) \quad \Leftrightarrow \quad 0\in \partial \varphi^{(a)}(\hat{u}).
\end{equation*} 
Given $\epsilon>0$, a point $\tilde{u}\in\realr^d$ is an $\epsilon-$approximation of $\hat{u}$ of type 1 whenever \cite[Definition 2.1]{Salzo2012}
\begin{equation}\label{eq:type1}
\tilde{u}\approx_{1}^{\epsilon} \hat{u} \quad \Leftrightarrow \quad 0\in\partial_{\epsilon} \varphi^{(a)}(\tilde{u}).
\end{equation} 
Likewise, a point $\tilde{u}\in\realr^d$ is an $\epsilon-$approximation of $\hat{u}$ of type 2 whenever \cite[Definition 2.2]{Salzo2012}
\begin{equation}\label{eq:type2}
\tilde{u}\approx_{2}^{\epsilon} \hat{u} \quad \Leftrightarrow \quad a-\tilde{u}\in \partial_{\epsilon} F(\tilde{u}).
\end{equation}
By using a well-known $\epsilon-$subdifferential calculus rule \cite[p. 1171]{Salzo2012}, we have that condition \eqref{eq:type2} implies condition \eqref{eq:type1}, which makes the notion of $\epsilon-$approximation of type 2 stronger than the previous one. For each $n\geq 0$, we define
\begin{align*}
\varphi_{\lambda,\mu_n}^{(n)}(u)&=\frac{1}{2}\|u-(u_n-\alpha\grad f(u_n)-\alpha \mu_n A^Tv_n)\|_2^2+\alpha \lambda g(u)\\
\varphi_{\lambda_n,\mu_n}^{(n)}(u)&=\frac{1}{2}\|u-(u_n-\alpha\grad f(u_n)-\alpha \mu_nA^Tv_n)\|_2^2+\alpha \lambda_n g(u)
\end{align*}
and
\begin{align*}
u_{\lambda,\mu_n}^{(n)}&=\underset{u\in\realr^d}{\operatorname{argmin}} \ \varphi_{\lambda,\mu_n}^{(n)}(u)=\operatorname{prox}_{\alpha\lambda g}(u_n-\alpha \grad f(u_n)-\alpha \mu_n A^Tv_n )\\
&=\operatorname{prox}_{\alpha\lambda g}(u_n-\alpha\grad f(u_n)-\alpha\mu A^Tv_n-\alpha(\mu_n-\mu)A^Tv_n ).
\end{align*}
We note that $u_{\lambda,\mu_n}^{(n)}$ is similar to the primal point $u_{n+1}$ defined in our method (2), except for the parameter $\lambda$ that replaces $\lambda_n$, and the error term $-\alpha(\mu_n-\mu)A^Tv_n$ inside the argument of the proximal operator. Then, we can write down the following implications
\begin{align*}
& u_{n+1} = \operatorname{prox}_{\alpha \lambda_n g}(u_n-\alpha\grad f(u_n)-\alpha\mu_nA^Tv_n) \ \Leftrightarrow \ u_{n+1} = \underset{u\in\realr^d}{\operatorname{argmin}} \ \varphi_{\lambda_n,\mu_n}^{(n)}(u)\\
&\Leftrightarrow \ 0\in\partial \varphi^{(n)}_{\lambda_n,\mu_n}(u_{n+1})\\
&\Leftrightarrow \ \varphi^{(n)}_{\lambda_n,\mu_n}(u)\geq \varphi^{(n)}_{\lambda_n,\mu_n}(u_{n+1}), \quad \forall \ u\in\realr^d\\
&\Rightarrow \ \varphi^{(n)}_{\lambda_n,\mu_n}(u_{\lambda,\mu_n}^{(n)})\geq \varphi^{(n)}_{\lambda_n,\mu_n}(u_{n+1})\\
&\Rightarrow \  \varphi^{(n)}_{\lambda,\mu_n}(u_{\lambda,\mu_n}^{(n)}) \geq \varphi^{(n)}_{\lambda,\mu_n}(u_{n+1})+\alpha(\lambda_n-\lambda)(g(u_{n+1})-g(u_{\lambda,\mu_n}^{(n)}))\\
&\Rightarrow \  \varphi^{(n)}_{\lambda,\mu_n}(u) \geq \varphi^{(n)}_{\lambda,\mu_n}(u_{n+1})-\alpha|\lambda_n-\lambda|\, |g(u_{n+1})-g(u_{\lambda,\mu_n}^{(n)})|, \quad \forall \ u\in \realr^d,
\end{align*}
where the last inequality follows from the fact that $u_{\lambda,\mu_n}^{(n)}$ is the unique minimum point of $\varphi_{\lambda,\mu_n}^{(n)}$. Since the operator  $T(\alpha,\lambda,\mu,u,v)=\operatorname{prox}_{\alpha \lambda g}(u-\alpha \grad f(u)-\alpha\mu A^Tv)$ is continuous with respect to $\mu$, $u$, $v$, the boundedness of $\{u_n\}_{n\in\naturaln}$ and $\{v_n\}_{n\in\naturaln}$ (guaranteed by inequality (25)) and the limit $\mu_n\rightarrow \mu$ imply that the sequence $\{u_{\lambda,\mu_n}^{(n)}\}_{n\in\naturaln}$ is bounded. Since $g$ is continuous by assumption, it follows that also $\{|g(u_{n+1})-g(u_{\lambda,\mu_n}^{(n)})|\}_{n\in\naturaln}$ is bounded. By letting $M_1=\sup_{n\in\naturaln}|g(u_{n+1})-g(u_{\lambda,\mu_n}^{(n)})|<\infty$, we have
\begin{align*}
u_{n+1} &= \operatorname{prox}_{\alpha \lambda_n g}(u_n-\alpha\grad f(u_n)-\alpha\mu_nA^Tv_n)\\
&\Rightarrow \ \varphi^{(n)}_{\lambda,\mu_n}(u) \geq \varphi^{(n)}_{\lambda,\mu_n}(u_{n+1})-\alpha M_1|\lambda_n-\lambda|, \quad \forall \ u\in \realr^d\\
& \Rightarrow \ 0\in\partial_{\epsilon_{n+1}}\varphi_{\lambda,\mu_n}^{(n)}(u_{n+1}), \quad \text{where }\epsilon_{n+1} =\alpha  M_1|\lambda_n-\lambda|\\
&\Rightarrow \ u_{n+1}\approx_1^{\epsilon_{n+1}}\operatorname{prox}_{\alpha \lambda g}(u_n-\alpha\grad f(u_n)-\alpha\mu_nA^Tv_n).
\end{align*}
Likewise, from the definition (\ref{eq:algorithm}) of the point $v_{n+1}$ and the corresponding optimality condition, we have
\begin{align*}
& v_{n+1} = \operatorname{prox}_{(\beta/\mu_n)h^*}(v_n+(\beta/\mu_n)A(2u_{n+1}-u_n))\\
& \Leftrightarrow \ v_n+(\beta\mu_n^{-1})A(2u_{n+1}-u_n)-v_{n+1} \in \partial (\beta\mu_n^{-1}h^*)(v_{n+1}).
\end{align*}
By applying the definition of subdifferential (4) to the previous optimality condition, we obtain the inequality
\begin{equation*}
(\beta\mu_n^{-1})h^*(v)\geq (\beta\mu_n^{-1})h^*(v_{n+1})+\langle v_n+(\beta\mu_n^{-1})A(2u_{n+1}-u_n)-v_{n+1},v-v_{n+1}\rangle, \quad 
\end{equation*}
$\forall \ v\in\operatorname{dom}(h^*)$.
By rearranging the terms involving $h^\ast$, we get for all $v\in\operatorname{dom}(h^*)$
\begin{align}
(\beta\mu^{-1})h^*(v)&\geq (\beta\mu^{-1})h^*(v_{n+1})+\langle v_n+(\beta\mu_n^{-1})A(2u_{n+1}-u_n)-v_{n+1},v-v_{n+1}\rangle\nonumber\\
& \quad+ \beta(\mu_{n}^{-1}-\mu^{-1})(h^*(v_{n+1})-h^*(v))\nonumber\\
& =  (\beta\mu^{-1})h^*(v_{n+1})+\langle v_n+(\beta\mu^{-1})A(2u_{n+1}-u_n)-v_{n+1},v-v_{n+1}\rangle\nonumber\\
& \quad + \beta(\mu_n^{-1}-\mu^{-1})(h^*(v_{n+1})-h^*(v) + \langle A(2u_{n+1}-u_n),v-v_{n+1}\rangle ),\nonumber\\
& \geq  (\beta\mu^{-1})h^*(v_{n+1})+\langle v_n+(\beta\mu^{-1})A(2u_{n+1}-u_n)-v_{n+1},v-v_{n+1}\rangle\nonumber\\
& \quad - \beta|\mu_n^{-1}-\mu^{-1}|(|h^*(v_{n+1})-h^*(v)| + \|A\|\, \|2u_{n+1}-u_n\|_2\, \|v-v_{n+1}\|_2 ),\label{eq:ine_dual}
\end{align}
where the last inequality is due to the application of the Cauchy-Schwarz inequality. Let us define 

\begin{equation*}
\mathcal{C}_{h^*,A}=\{|h^*(v_{n+1})-h^*(v)| + \|A\|\,\|2u_{n+1}-u_n\|_2\,\|v-v_{n+1}\|_2: \ n\in\naturaln, v\in\operatorname{dom}(h^*)\}\subseteq \realr_{\geq 0}.
\end{equation*}
Since $\{u_n\}_{n\in\naturaln}$, $\{v_n\}_{n\in\naturaln}$, $\operatorname{dom}(h^*)$ are bounded and $h^*$ is continuous on its domain, we conclude that the set $\mathcal{C}_{h^*,A}$ is bounded. By letting $M_2 = \sup \mathcal{C}_{h^*,A} < \infty$, inequality \eqref{eq:ine_dual} implies that
\begin{align*}
(\beta\mu^{-1})h^*(v)&\geq (\beta\mu^{-1})h^*(v_{n+1})+\langle v_n+(\beta\mu^{-1})A(2u_{n+1}-u_n)-v_{n+1},v-v_{n+1}\rangle\nonumber\\
&-\beta M_2|\mu_n^{-1}-\mu^{-1}|, \quad \forall \ v\in\operatorname{dom}(h^*).
\end{align*}
By definition \eqref{eq:eps_sub} of $\epsilon-$subdifferential, the previous inequality yields
\begin{align*}
v_n+(\beta\mu^{-1})A(2u_{n+1}-u_n)-v_{n+1} \in \partial_{\delta_{n+1}} (\beta\mu^{-1}h^*)(v_{n+1}),
\end{align*}
where $\delta_{n+1} = \beta M_2|\mu_n^{-1}-\mu^{-1}|$.
In accordance with the definition of $\epsilon-$approximation of type 2 given in \eqref{eq:type2}, the above subdifferential inclusion is equivalent to writing
\begin{equation*}
v_{n+1}\approx_2^{\delta_{n+1}} \operatorname{prox}_{(\beta/\mu)h^*}(v_n+(\beta/\mu)A(2u_{n+1}-u_n)).
\end{equation*}
Therefore, under the aforementioned assumptions on $g$ and $h^*$, we have shown that method (2) can be reformulated as 
\begin{align*}
e_{n+1}&=(\mu_n-\mu)A^Tv_n\\
u_{n+1}&\approx_1^{\epsilon_{n+1}}\operatorname{prox}_{\alpha \lambda g}(u_n-\alpha\grad f(u_n)-\alpha\mu A^Tv_n-\alpha e_{n+1}), \\[2mm]
&\qquad\qquad\qquad\qquad\qquad\qquad\qquad\qquad \text{where }\epsilon_{n+1}=\alpha M_1|\lambda_n-\lambda|\\
v_{n+1}&\approx_2^{\delta_{n+1}} \operatorname{prox}_{(\beta/\mu)h^*}(v_n+(\beta/\mu)A(2u_{n+1}-u_n)), \\[2mm]
&\qquad\qquad\qquad\qquad\qquad\qquad\qquad\qquad \text{where }\delta_{n+1}=\beta M_2|\mu_n^{-1}-\mu^{-1}|.
\end{align*}
As it is written above, method (2) is a special instance of the inexact primal--dual method proposed in \cite[Section 3.1]{Rasch2020}, corresponding to the case where the method is applied to the saddle point problem $\min_{u\in\realr^d}\max_{v\in\realr^{d'}}f(u)+\lambda g(u) +\langle \mu Au,v\rangle - \mu h^*(v)$ with a specific choice of the parameters $e_{n+1}$, $\epsilon_{n+1}$, $\delta_{n+1}$. Consequently, the convergence analysis in \cite[Section 3.1]{Rasch2020} is applicable to our method, however under more restrictive assumptions than the ones required in Theorem 3.5. Indeed, according to \cite[Theorem 2]{Rasch2020}, one could prove the convergence of method (2) under a tighter bound on $\alpha$ than the one proposed in Theorem 3.5, and assuming that $\sum_{n\in\naturaln}\|e_{n}\|<\infty$, $\sum_{n\in\naturaln}\sqrt{\epsilon_n}<\infty$ and $\sum_{n\in\naturaln}\delta_n<\infty$; the first of these conditions requires that $\{v_n\}_{n\in\naturaln}$ is bounded (which is proved in Theorem 3.5) and $\sum_{n\in\naturaln}|\mu_n-\mu|<\infty$, whereas the second one is ensured only if $\sum_{n\in\naturaln}\sqrt{|\lambda-\lambda_n|}<\infty$, which is a stronger assumption than the one made in Theorem 3.5. By contrast, our result is applicable to more general problems under less restrictive assumptions on the parameters.

\section{Numerical illustration}\label{sec:numerics}

\subsection{Total variation deblurring}

\begin{figure}
\centering\includegraphics{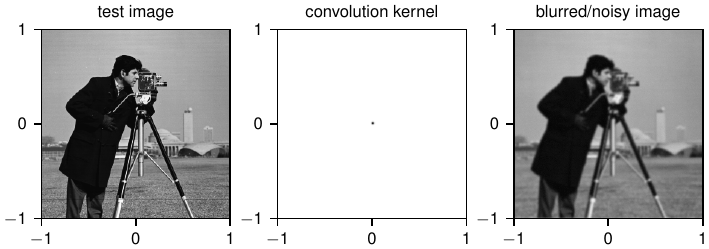}
\caption{From left to right: ground truth, convolution kernel and noisy image.}\label{fig:pics}
\end{figure}

In this Section we present a numerical optimization experiment in the context of Total Variation image reconstruction \cite{Chambolle2016,Rudin1992}. We consider a ground truth image (Figure~\ref{fig:pics}) which is convolved with a Gaussian kernel and contaminated by Gaussian noise to obtain synthetic blurred and noisy data $y$:
\begin{displaymath}
y=K u^\text{input}+\text{noise}.
\end{displaymath}
\amendment{The variance of the Gaussian blur is $\sigma^2=10^{-4}$ and the level of the noise is $1\%$ of the noiseless data: $\|\text{noise}\|_2= \|K u^\text{input}\|_2/100$.}

The associated optimization problem we wish to solve, in order to recover an approximation of $u^\text{input}$ from $y$ and $K$, is
\begin{equation}\label{eq:tv-problem}
\min_{u\in\realr^d}\frac{1}{2}\|Ku-y\|_2^2+\mu \TV(u)+\ind_{[0,1]}(u).
\end{equation}
This is a special instance of problem (\ref{eq:problem}) with $f(u)=\frac{1}{2}\|Ku-y\|_2^2$, $h(v)=\|v\|_{1,2}$, $A$ is a local image gradient operator and $g=\ind_{[0,1]}$ constrains  each image pixel value between $0$ and $1$.

\begin{figure}
\centering\includegraphics{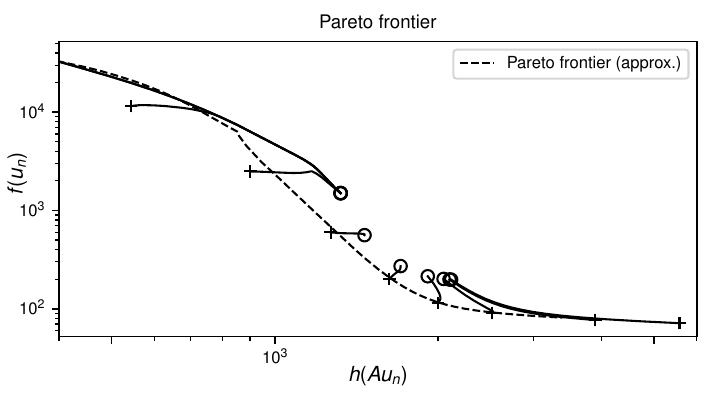}
\caption{Approximation of the Pareto-frontier by Algorithm (\ref{eq:algorithm}) in the context of the numerical experiment of Section~\ref{sec:numerics}. The dashed line represents the path followed by the iterates of proposed Algorithm (\ref{eq:algorithm}). The $10$ solid lines, on the other hand, correspond to the paths followed by the classical algorithm (\ref{eq:alg1}) for $10$ values of $\mu$ between $10^{3}$ and $10^{-3}$. The circles indicate the positions of the iterates (in the $f-h$ plane) near the start of the iterative process ($n=10$) and the crosses indicate the positions at end of the iterative process ($n=1000$). One sees that the dashed line is indeed a good approximation of the Pareto-frontier for a wide range of the penalty parameter $\mu$, as compared to its samples obtained from running algorithm (\ref{eq:alg1}) many times. A logarithmic scale is used on both axes.}\label{fig:pareto}
\end{figure}

In order to illustrate the usefulness of the optimization-by-continuation algorithm (\ref{eq:algorithm}) for computing the Pareto frontier for this problem, we first run the classical algorithm (\ref{eq:alg1}) for $10$ different values of the parameter $\mu$ (evenly spaced on a logarithmic scale from $10^3$ to $10^{-3}$). 
\amendment{The initial guess of the solution is always chosen as the zero image ($u_0=0$, $v_0=0$).}
The results of each run of $1000$ iterations are plotted in the trade-off plane, i.e. the values of $f(u_n)$ versus $h(u_n)$ are plotted as solid lines on Figure~\ref{fig:pareto} for each of the $10$ choices of $\mu$. Therefore, $10.000$ iterations are required to sample the Pareto frontier at $10$ points (end points of each solid line).

Secondly, we run the continuation method (\ref{eq:algorithm}) once, also for $1000$ steps, and starting from the previously obtained minimizer corresponding to $\mu=10^3$. In this case, the penalty parameters $\mu_n$ (on which algorithm (\ref{eq:algorithm}) depends) are also chosen between $10^3$ to $10^{-3}$ ($1000$ evenly spaced values, on a logarithmic scale). Hence this run requires $2000$ iterations in total (1000 for the initial point and 1000 for all other approximate minimizers). The result of this run is also plotted on Figure~\ref{fig:pareto} (dashed line). It is observed empirically that the dashed line offers a valid approximation of the Pareto-curve for the problem at hand.

\begin{figure}
\centering\includegraphics{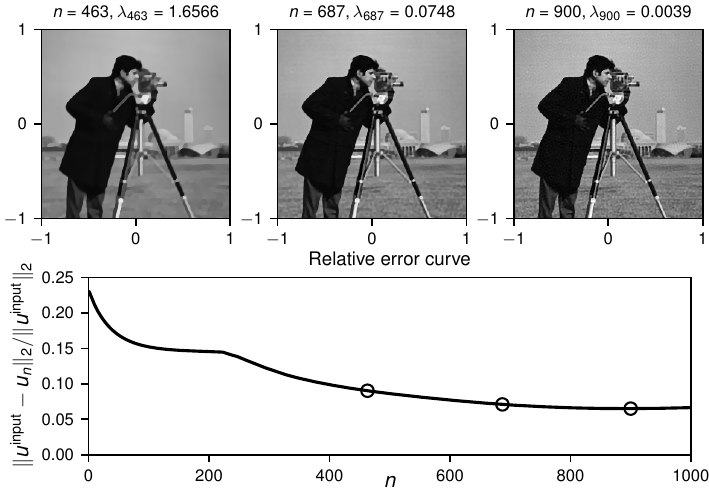}
\caption{\amendment{Top: Three intermediate iterates ($n= 463, 687, 900$) of Algorithm (\ref{eq:algorithm}) for the cameraman example. Bottom: the relative error $\|u^{\text{input}}-u_n\|_2/\|u^{\text{input}}\|_2$ as a function of $n$.}}
\label{fig:error}
\end{figure}

\amendment{
Figure~\ref{fig:error} depicts three \emph{intermediate} iterates of algorithm (\ref{eq:algorithm}) applied to the cameraman deblurring problem ($n= 463, 687, 900$, corresponding to $\lambda_{463}=1.66, \lambda_{687}=0.075, \lambda_{900}=0.0039$ respectively). The value of the corresponding data misfits $\|Ku_n-y\|_2/\|\text{noise}\|$ is $2$, $1$ and $0.9$ respectively. According to the Morozov discrepancy principle \cite{Morozov1966,Engl2000,Bonesky2009}, one should choose the penalty parameter $\lambda$ in problem (\ref{eq:tv-problem}) such that $\|K\hat u-y\|_2/\|\text{noise}\|\approx 1$ (with $\hat u$ the minimizer of problem (\ref{eq:tv-problem})). In this case, according to the discrepancy principle, we should prefer solution $u_{687}$. Indeed, by looking carefully at the first row of Figure~\ref{fig:error} one sees that image $u_{463}$ has too little detail (underfit). Image $u_{900}$ is also acceptable. In this particular example, the L-curve method \cite{Hansen2001} does not work very well as there is no distinctive bend in the Pareto-curve in Figure~\ref{fig:pareto}.
}

\amendment{
\subsection{Total Variation image deblurring with two penalty parameters}\label{sec:numerics2}

In this subsection we develop an example where two penalty parameters are present in the cost function. In the example of the previous subsection, there was effectively only one penalty parameter because of the use of an indicator function as the function $g$.

In this numerical example we apply the proposed algorithm to a different image deblurring problem. We now use a synthetic ground truth image (similar to the one in \cite{Jin2009}) which is again convolved with a Gaussian kernel and contaminated by Gaussian noise to obtain synthetic blurred and noisy data $y$:
\begin{displaymath}
y=K u^\text{input}+\text{noise};
\end{displaymath} 
see Figure~\ref{fig:synground}. The variance of the Gaussian blur is $\sigma^2=10^{-3}$ and the level of the noise is $10\%$ of the noiseless data: $\|\text{noise}\|_2= 0.10 \times \|K u^\text{input}\|_2$.

In this case, we illustrate the proposed algorithm on the optimization problem 
\begin{equation}\label{eq:tv-problem2}
\min_{u\in\realr^d}\frac{1}{2}\|Ku-y\|_2^2+\mu \TV(u)+\lambda \left(\|u\|_1+\ind_{[0,1]}(u)\right),
\end{equation}
i.e. a combination of a least squares data misfit term and \emph{two} penalty terms exactly as in problem (\ref{eq:problem}): $f(u)=\frac{1}{2}\|Ku-y\|_2^2$, $h(v)=\|v\|_{1,2}$, $A$ is a local image gradient operator and $g(u)=\|u\|_1+\ind_{[0,1]}(u)$ is a sparsity promoting $1$-norm (limited to the $[0,1]$ interval).
The joint use of Total Variation and $\ell_1$-norm regularization has been employed in several papers related to tomography, e.g. \cite{Dutta2012,Tang2017,Tong2018}.

\begin{figure}
\centering\includegraphics{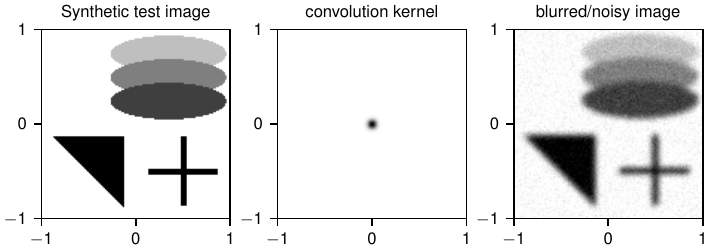}
\caption{\amendment{From left to right: a synthetic test image (see also \cite{Jin2009}), a Gaussian convolution kernel and the corresponding blurred and noisy test image.}}
\label{fig:synground}
\end{figure}

In order to show the efficacy of algorithm (\ref{eq:algorithm}) for computing (approximations) of several minimizers of cost function (\ref{eq:tv-problem2}) at once, we choose $20$ values of $\lambda$ and $\mu$ which are logarithmically spaced from $1$ to $0.02$ (for $\lambda$) and from $10$ to $1$ (for $\mu$). Only $20$ iterations of algorithm (\ref{eq:algorithm}) are performed, which is reasonable given the simplicity of the synthetic test image w.r.t the TV penalty. The starting point of the iteration is the zero image (which is the exact minimizer of (\ref{eq:tv-problem2}) for large $\lambda$ due to the sparsity promotion of the $1$-norm term). Three out of the twenty approximations are displayed in Figure~\ref{fig:contalg2}. As expected, large values of $\lambda$ lead to poor reconstruction results (the ground truth is not sparse).

\begin{figure}
\centering\includegraphics{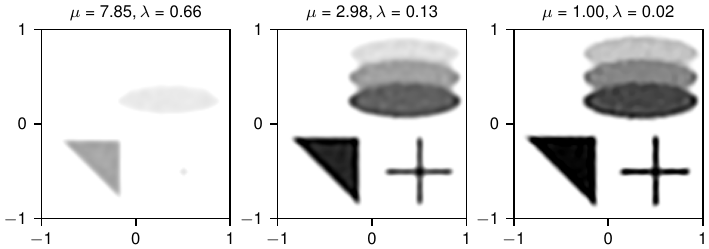}
\caption{\amendment{Three different reconstructions of the synthetic test image of Figure~\ref{fig:synground} for three different choices of the penalty parameters in model (\ref{eq:tv-problem2}). These reconstructions were obtained using algorithm (\ref{eq:algorithm}). The reconstruction on the left does not fit the data very well, whereas the one on the right is to be preferred based on the discrepancy principle: $\|K\hat u-y\|_2/\|\text{noise}\|\approx 1$.}}
\label{fig:contalg2}
\end{figure}

In Figure~\ref{fig:syntheticreconstruction} we show three reconstructions for the same choices of penalty parameters as in Figure~\ref{fig:contalg2} but using the traditional algorithm (\ref{eq:alg1}). Here we used $1000$ iterations starting each time from the zero image. At the price of more iterations, a subjective better image quality is obtained for the reconstruction on the right.
This is also confirmed on Figure~\ref{fig:error2} which shows  the relative error of each of the twenty reconstructions (for the twenty different values of the two penalty parameters) of algorithm (\ref{eq:algorithm}). It also shows the relative error for the three reconstructions of Figure~\ref{fig:syntheticreconstruction} obtained from algorithm (\ref{eq:alg1}) after $1000$ iterations.

Comparing $20$ iterations of algorithm (\ref{eq:algorithm}) with $1000$ iterations of algorithm (\ref{eq:alg1}), the numerical improvement provided by algorithm (\ref{eq:alg1}) on the relative error is modest. Hence, we conclude that for the price of $20$ iterations, the algorithm (\ref{eq:algorithm}) can approximate well $20$ minimizers of problem (\ref{eq:tv-problem2}) for the $20$ different values of the two penalty parameters. Running $20$ instances of  algorithm (\ref{eq:alg1}) would be more expensive, even if each instance runs for just $10$ iterations. On the other hand, once the ``good'' values of the penalty parameters are known, the traditional algorithm (\ref{eq:alg1}) can be used to refine the corresponding solution.

It is not claimed that the algorithm (\ref{eq:algorithm}) is faster than algorithm (\ref{eq:alg1}) if the goal is to obtain the minimizer of (\ref{eq:tv-problem2}) for a single fixed value of the penalty parameters; this kind of speed-up is not our aim here (see e.g. \cite{Hale2008}), rather we speed up the computation of a whole set of minimizers at once.

%l1 norm deblurring is similar to halftone as in newspaper pictures (gray scale imitation with black/white...

\begin{figure}
\centering\includegraphics{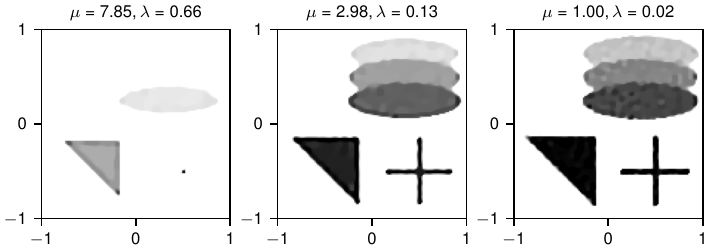}
\caption{\amendment{Three reconstructions of the synthetic image of Figure~\ref{fig:synground} for the same three choices of the two penalty parameters as in Figure~\ref{fig:contalg2}, now obtained with $1000$ iterations of algorithm (\ref{eq:alg1}) each.}}
\label{fig:syntheticreconstruction}
\end{figure}

\begin{figure}
\centering\includegraphics{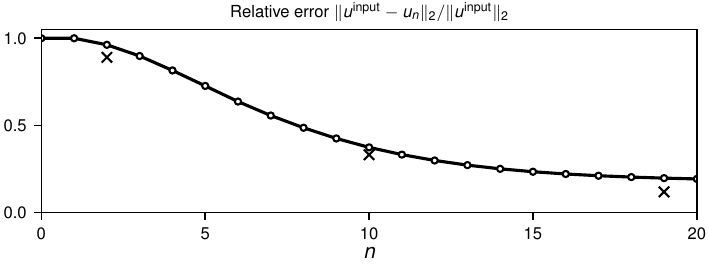}
\caption{\amendment{The relative error $\|u^{\text{input}}-u_n\|_2/\|u^{\text{input}}\|_2$ for the $20$ different reconstructions $u_n$ obtained from algorithm (\ref{eq:algorithm}). Crosses indicate the relative error of the three reconstructions of Figure~\ref{fig:syntheticreconstruction} obtained after $1000$ iterations of algorithm (\ref{eq:alg1}). This Figure illustrates the use of the algorithm (\ref{eq:algorithm}) for computing a reasonable approximation of the minimizers of (\ref{eq:tv-problem2}) for a whole set of penalty parameters at once.}}
\label{fig:error2}
\end{figure}

}%end amendment

\section{Conclusions}\label{sec-concl}

An optimization-by-continuation algorithm is proposed for problem (\ref{eq:problem}). It is a first order primal-dual method with variable penalty parameters (not variable step lengths). A proof of convergence employing a variable metric criterion was given. Such methods are usually proposed for discussing convergence of Newton-like accelerated first-order optimization algorithms.

Under some additional assumptions on the cost function, the proposed algorithm can be interpreted as a special instance of an inexact proximal algorithm proposed in \cite{Rasch2020}. It follows that the convergence analysis in \cite{Rasch2020} is also applicable to the proposed algorithm, however under more restrictive assumptions.

Future work, not addressed in this paper, deals with the convergence of the generalized ISTA algorithm \cite{Loris2011,Chen2019} with variable penalty parameters. Another line of research may focus on extending the results here to the more general class of forward backward splitting algorithms, which are suitable for problems formulated in terms of monotone operators without reference to optimization problems \cite{Bauschke2011}.

Finally, another topic, also not addressed in this paper, concerns the behavior of the iterates when $N$ penalty parameters $\lambda_n$ ($n=1,\ldots, N$) are chosen (e.g. uniformly and decreasing) in a fixed interval $[\lambda_{\min},\lambda_{\max}]$ as is the case in practice. What can be said about the iterates when $N$ tends to infinity? It would indeed be useful to prove that the set of iterates $\{u_n^{(N)}, n\in [0,N]\}$ approaches the Pareto frontier when $N\to\infty$, or to quantify the distance between those two sets.

\section*{Acknowledgments}

This work was supported by the Fonds de la Recherche Scientifique - FNRS under Grant CDR J.0122.21.
\amendment{Simone Rebegoldi is a member of the INDAM research group GNCS and is supported by the Italian MUR through the PRIN 2022 project ``Inverse problems in PDE: theoretical and numerical analysis'', project code: 2022B32J5C (CUP B53D23009200006), and the PRIN 2022 PNRR Project ``Advanced optimization METhods for automated central veIn Sign detection in multiple sclerosis from magneTic resonAnce imaging (AMETISTA)'', project code: P2022J9SNP (CUP  E53D23017980001), under the National Recovery and Resilience Plan (PNRR), Italy, Mission 04 Component 2 Investment 1.1 funded by the European Commission - NextGeneration EU programme.}
%old:
%Simone Rebegoldi is a member of the INDAM research group GNCS and is partially supported by the PRIN project 2022B32J5C, funded by the Italian Ministry of University and Research.

%\bibliography{ParetoPD}{}
%\bibliographystyle{elsarticle-num.bst}

\end{document}